\documentclass[a4paper, 10pt]{article}
 
\usepackage{enumerate}
\usepackage{authblk}
\usepackage[latin1]{inputenc}
\usepackage[T1]{fontenc}
\usepackage{amsfonts}
\usepackage{graphicx}
\usepackage{ulem}
\usepackage{lmodern}
\usepackage[top=4cm, bottom=4cm, left=4cm, right=4cm]{geometry}
\usepackage{amsthm}
\usepackage{amsmath}
\usepackage{amssymb}
\usepackage{mathrsfs}
\usepackage{enumitem}
\usepackage{sectsty}
\usepackage{amssymb, url}
\usepackage[colorlinks=true]{hyperref} 
\usepackage{verbatim}
\usepackage{fancyhdr}
\allsectionsfont{\centering}

\newenvironment{ack}{\bigskip \noindent {\bf Acknowledgement}\it}

\theoremstyle{definition}
\newtheorem{Theorem}{Theorem}
\newtheorem*{Theorem*}{Theorem}
\newtheorem{Theoremroman}{Theorem}

\newtheorem*{Remark*}{Remark}
\newtheorem{Lemma}[Theorem]{Lemma}

\newtheorem{Proposition}[Theorem]{Proposition}
\newtheorem*{Proposition*}{Proposition}

\theoremstyle{remark}
\newtheorem{Remark}[Theorem]{Remark}

\newcommand{\ord}{\operatorname{ord}}

\newcommand{\id}{\operatorname{id}}

\newcommand{\Gal}{\operatorname{Gal}}

\newcommand{\Mod}[1]{\ (\mathrm{mod}\ #1)}

\newcommand{\SLZ}{\operatorname{SL}_2(\bZ)}
\newcommand{\Span}{\operatorname{Span}}
\newcommand\smallmat[4]{\left(\begin{smallmatrix}#1&#2\\#3&#4\end{smallmatrix}\right)}

\newcommand{\bC}{\mathbb{C}}

\newcommand{\bP}{\mathbb{P}}
\newcommand{\bQ}{\mathbb{Q}}
\newcommand{\bR}{\mathbb{R}}

\newcommand{\bZ}{\mathbb{Z}}

\newcommand{\cE}{\mathcal{E}}

\newcommand{\cH}{\mathcal{H}}

\newcommand{\cM}{\mathcal{M}}

\newcommand{\cO}{\mathcal{O}}

\newcommand{\cS}{\mathcal{S}}

\newcommand{\fc}{\mathfrak{c}}

\begin{document}

\title{\textbf{Generalized modular forms with a cuspidal divisor}}
\author{{Quentin Gazda}\thanks{Current address: 
Univ Lyon, Universit\'e Jean Monnet Saint-\'Etienne, CNRS UMR 5208, Institut Camille Jordan, F-42023 Saint-\'Etienne, France}}
\date{June, 2020}

\maketitle

\begin{abstract}
In \cite{WinKoh}, Kohnen proves that if $\Gamma=\Gamma_0(N)$ where $N$ is a square-free integer, then any modular function of weight $0$ for $\Gamma$ having a divisor supported at the cusps is an $\eta$-quotient. Under the condition of having rational Fourier coefficients, we are able to extend Kohnen's result to the case where $N$ is the square of a prime. If the rationality condition does not hold, we show that the statement is no longer true by providing a family of counter-examples that generalizes naturally the Dedekind $\eta$-function. This paper fits within the framework of generalized modular forms in the sense of Knopp and Mason. 
\end{abstract}

\section{Introduction}
The theory of generalized modular forms was initiated by Knopp and Mason in \cite{KnoMas}. It consists in holomorphic functions on $\cH$, meromorphic at the cusps, which satisfies a modular invariance property for multiplier systems not necessarily of finite order. The main motivation for introducing them arises from \textit{Moonshine} properties for finite groups in \textit{orbifolds theory}, where generalized modular forms naturally appear as trace elements. This is discussed in \cite{KohMas}. From the latter reference, a main concern is how to detect classical modular forms out of generalized one. As in the classical case, generalized modular forms can be expanded in Fourier series. It was conjectured by Kohnen and Mason that normalized generalized modular form with integral Fourier coefficients are classical (see \cite[(2)]{KohMas}). \\

The goal of this paper is to discuss various particular cases of this conjecture when generalized modular forms have a divisor supported at the cusps of $\Gamma=\Gamma_0(N)$. To state our result, a brief preliminary discussion on notation is called for. \\

Let $g$ be any function from the Poincar\'e upper half-plane $\cH$ to $\bC$, $k$ an integer and $\sigma=\smallmat{a}{b}{c}{d}$ in $\operatorname{GL}_2(\bQ)$ with $ad-bc>0$. We denote by $(g|_k \sigma)$ the weight $k$ slash operator of $\sigma$, defined by:
\begin{equation}
(g|_k\sigma)(\tau)\stackrel{\text{def}}{=}(c\tau+d)^{-k}g\left(\frac{a\tau+b}{c\tau+d}\right)\quad (\tau\in \cH). \nonumber
\end{equation}
Let $\Gamma$ be a finite index subgroup of $\SLZ$. A \textit{generalized modular function (of weight $0$) with a cuspidal divisor $f$ for} $\Gamma$ is a meromorphic function $f$ on $\cH$ such that
\begin{enumerate}[label=$(\alph*)$]
\item As a meromorphic function, $f$ has an empty divisor on $\cH$,
\item \label{item:b} For every $\sigma \in \Gamma$, the quotient $(f|_0\sigma)/f$ is constant as a function on $\cH$,
\item \label{item:c} at each cusps $\fc=a/c$ in $\bP^1(\bQ)$ written in lowest term ($c$ might be zero in which case $\fc=\infty$), $f$ has an expansion of the form
\begin{equation}\label{expansion-at-cusps}
(f|_0\gamma)(\tau)=\sum_{n\geq n_{\fc}(f)}{a_{\gamma,n}(f)q^{n/w_{\fc}}} \quad (q=e^{2i\pi \tau},~\tau\in \cH)
\end{equation}
for $\gamma\in \SLZ$ having $a$ and $c$ in its first column, where $n_{\fc}(f)$ is an integer with $a_{\gamma,n_{\fc}}(f)\neq 0$, and $w_{\fc}$ is the width of $\fc$ for $\Gamma$, namely, the smallest positive integer $w$ such that $\smallmat{1}{w}{0}{1}$ belongs to $\gamma^{-1}\Gamma \gamma$. 
\end{enumerate}
In \cite{KnoMas}, an holomorphic function $f$ on $\cH$ satisfying \ref{item:b} and \ref{item:c} is rather called a \textit{parabolic generalized modular form}. Note that the only condition that extends the definition of modular functions (of weight $0$ for $\Gamma$ with a cuspidal divisor) is condition \ref{item:b}.\\

As in the case of modular functions, the complex number $a_n(f):=a_{\id,n}(f)$ will be referred to as the $n^{th}$-\textit{Fourier coefficient} of $f$. Also, given a cusp $\fc$, we define the order $\ord_{\fc}(f)$ of $f$ at $\fc$ to be the integer $n_{\fc}(f)$ appearing in formula \eqref{expansion-at-cusps}. If the first coefficient of $f$ -- namely the $\ord_{\infty}(f)^{th}$-Fourier coefficient -- is one, we say that $f$ is \textit{normalized}. \\
To any $f$, we associate a character $\nu_f:\Gamma\to \bC^{\times}$ of $\Gamma$ by $\nu_f(\gamma)=(f|_0\gamma)/f$. We call $\nu_f$ the \textit{multiplier system of $f$}. In general, it might not be of finite order, nor unitary (i.e. $|\nu_f(\gamma)|=1$, $\forall \gamma \in \Gamma$).\\

Here and from now on $\tau\in \cH$. It was prove by Kohnen in \cite[Thm.~2]{WinKoh} that if $\Gamma=\Gamma_0(N)$ where $N$ is a square-free integer and $f$ has a multiplier system of finite order, then $f$ is an $\eta$-quotient, i.e. there exists $\alpha$ in $\bC$ and rational numbers $(c_t)_{t|N}$ such that 
\begin{equation}
f(\tau)=\alpha \prod_{t|N}{\eta(t\tau)^{c_t}}, \quad \eta(\tau)\stackrel{\text{def}}{=} q^{\frac{1}{24}}\prod_{n=1}^{\infty}{(1-q^n)}, \nonumber
\end{equation} 
where $\eta$ is the {\it Dedekind eta function} and where rational powers are defined in terms of the principal branch of the complex logarithm. It is stated in \cite{WinKoh}:
\begin{quote}
\textit{We think that the assertion of Theorem $2$ or a slightly weaker statement eventually would be true for arbitrary $N$ $[...]$.}
\end{quote}
Under a rationality condition, we are able to extend Kohnen's result in the case where $N$ is the square of a prime number. We prove:
\begin{Theoremroman}\label{thmI}
Let $p$ be a prime number and let $f$ be a normalized generalized modular function with a cuspidal divisor for $\Gamma_0(p^2)$ having rational coefficients. The following are equivalent:
\begin{enumerate}[label=$(\roman*)$]
\item \label{item:i}$f$ is an $\eta$-quotient, i.e there exist rational numbers $a$, $b$ and $c$ such that 
\begin{equation}
f(\tau)=\eta(\tau)^{a}\eta(p\tau)^{b}\eta(p^2\tau)^c, \nonumber
\end{equation}
\item \label{item:ii} $f$ has integral Fourier coefficients,
\item \label{item:iii} $f$ has $\ell$-integral Fourier coefficients for almost all prime $\ell$,
\item \label{item:iv} the multiplier system $\nu_f$ has finite order,
\item \label{item:v} the multiplier system $\nu_f$ is unitary.
\end{enumerate}
\end{Theoremroman}
\begin{Remark}
Note that Theorem \ref{thmI} is easily extended to weight $k$ by multiplying with the correct power of $\eta$.
\end{Remark}
We can achieve the proof of the above theorem thanks to an explicit computation, but the general case when $p^2$ is replaced by $N$, although likely to be true, would follow from a different argument. \\

It also appears that if the rationality condition is not fulfilled, the previous does not hold anymore. As a counterexample, we provide:
\begin{Theoremroman}\label{thmII}
The following product defines a generalized modular function with a cuspidal divisor for $\Gamma_0(25)$:
\begin{equation}
\eta_{\left(\frac{\cdot}{5}\right)}(\tau)=\prod_{n=1}^{\infty}{\left(\frac{1-\frac{1-\sqrt{5}}{2}q^n+q^{2n}}{1-\frac{1+\sqrt{5}}{2}q^n+q^{2n}}\right)^{\left(\frac{n}{5}\right)}}
\end{equation}
where $\left(\frac{\cdot}{5}\right)$ denotes the Legendre symbol for $5$. Its multiplier system has finite order and $\eta_{\left(\frac{\cdot}{5}\right)}$ is not an $\eta$-quotient. 
\end{Theoremroman}
In the sequel, generalizations $\eta_{\chi}$ of the previous function to any real and primitive Dirichlet character $\chi$ is introduced (see Proposition \ref{holo}). We leave as an open question whether its multiplier system $\nu_{\eta_{\chi}}$ is of finite order or not in general. \\

The main technic in the proof of the two above theorems relies on \textit{Knopp and Mason's Theorem} \cite[Thm.~2]{KnoMas}. Let $G(\Gamma)$ be the set of generalized modular functions with a cuspidal divisor for $\Gamma$. This is an abelian group with the product of functions. The logarithmic differentiation $\theta$ defines a morphism of abelian groups:
\begin{equation}
\theta:G(\Gamma)\longrightarrow \cM_2(\Gamma),\quad f\longmapsto \frac{1}{2i\pi}\frac{df}{f},
\end{equation}
where $\cM_2(\Gamma)$ is the space of modular forms of weight $2$ for $\Gamma$ and where $d$ denotes the derivation in the variable $\tau$. In \cite[Thm.~2]{KnoMas}\footnote{Remember that what we call \textit{generalized modular form} is qualified \textit{parabolic} in \cite{KnoMas}.}, it is inferred that the kernel of $\theta$ is $\bC^{\times}$ and that its image consists in elements of $\cM_2(\Gamma)$ whose constant term at any \textit{finite} cusp $\fc=a/c$ (where $c\neq 0$) belongs to $(w_{\fc}c^2)^{-1}\bZ$. In more details, an element $f$ as in equation \eqref{expansion-at-cusps} will be sent under $\theta$ to an element in $\cM_2(\Gamma)$ having constant term at $\fc$ given by $(w_{\fc}c^2)^{-1}\ord_{\fc}(f)$. Our main results will follow from a careful computation of the constant terms at each cusp of $\Gamma=\Gamma_0(N)$ of the Hecke eigenbasis of the subspace of Eisenstein series $\cE_2(\Gamma)$. \\

The paper is organized as follow. In Section \ref{sec:preliminaries} we present the Hecke eigenbasis of the space of weight $2$ Eisenstein series according to \cite[Section~4.6]{Diamond}. In Section \ref{sec:results} we use the basis presented in Section \ref{sec:preliminaries} together with the map $\theta$ to prove Theorem \ref{extension} which extends slightly Kohnen's \cite[Thm.~2]{WinKoh}. Proposition \ref{expansion-at-cusp} will be of central importance for our work regarding Knopp and Mason's \cite[Thm.~2]{KnoMas} as it computes the constant term at the cusps of the previous Eisenstein eigenforms. The rest of the paper consists in the proof of the above theorems.

\begin{ack}
In the preparation of an early version of this work, the author benefits from multiple discussions with Winfried Kohnen. May he find here the expression of my profound gratitude. I am also much obliged to the referee for his careful reading.
\end{ack}

\section{Preliminaries}\label{sec:preliminaries}
Let $\Gamma$ be a congruence subgroup of $\operatorname{SL}_2(\bZ)$. We start by listing some subgroups of $G(\Gamma)$ of particular interest. 
\begin{enumerate}
\item The subgroup of $G(\Gamma)$ consisting of $f$ for which $\ord_{\fc}(f)=0$ for each cusp $\fc$ of $\Gamma$ will be denoted $G_0(\Gamma)$.
\item The subgroup of $G(\Gamma)$ consisting of $f$ for which $\nu_f$ has finite order will be denoted $G_c(\Gamma)$.  The subgroup consisting of $f$ for which $\nu_f$ is unitary will be denoted $G_b(\Gamma)$. The subscripts $c$ and $b$ stand respectively for \textit{classical} and \textit{bounded}.
\item The subgroup of elements that are normalized and have integral coefficients will be denoted $G_{\bZ}(\Gamma)$.
\end{enumerate}
We warn the reader that the inclusion $\bC^{\times}\subset G_0(\Gamma)$ might not be an equality. Any $f\in G_0(\Gamma)$ is indeed holomorphic on $\cH$ and at all cusps but, as $|\nu_f(\Gamma)|\subset \bR$ might not be bounded, one cannot conclude that $|f|$ is bounded on $\cH$. What however is true from Liouville's Theorem is the equality
\begin{equation} \label{liouville}
G_0(\Gamma)\cap G_b(\Gamma)=G_0(\Gamma)\cap G_c(\Gamma)=\bC^{\times}
\end{equation}
where $G_c(\Gamma)\subset G_b(\Gamma)$ might also not be an equality in general. From \cite[Thm.~2]{KnoMas}, its proof and its supplement, we in addition have an exact sequence of abelian groups:
\begin{equation}
0\longrightarrow \bC^{\times} \longrightarrow G_0(\Gamma) \stackrel{\theta}{\longrightarrow} \cS_2(\Gamma)\longrightarrow 0. \nonumber
\end{equation}
Any $\tau_0\in \cH$ defines a section $s:\cS_2(\Gamma)\to G_0(\Gamma)$ mapping $f$ to 
\begin{equation}
s(f)(\tau)=\exp\left(2i\pi \int_{\tau_0}^{\tau}{f(u)du}\right) \nonumber
\end{equation}
where we integrate over any path in $\cH$ from $\tau_0$ to $\tau$. \\
It is also clear that $\theta(G_{\bZ}(\Gamma))\subset\cM_2(\Gamma,\bZ)$, the latter being the group of modular forms of weight $2$ having integral Fourier coefficients. Kohnen and Mason's conjecture \cite[(2)]{KohMas} should also imply $G_{\bZ}(\Gamma)\subset G_c(\Gamma)$ which remains an open question. \\

Let $N$ be a positive integer and let $\Gamma=\Gamma_0(N)$. For any $t|N$, let $V_t$ be $\smallmat{t}{0}{0}{1}$. The function $(\eta|_0V_t)/\eta$ is an element of $G_c(\Gamma)\cap G_{\bZ}(\Gamma)$ (see \cite[Prop.~1.2.5]{Diamond}). We denote $G_{\eta}(\Gamma)$ the subgroup of $G(\Gamma)$ generated by those elements. Generally, the latter is referred to as the {\it group of $\eta$-quotients}. Let
\begin{equation}
G_2(\tau)\stackrel{\text{def}}{=}\frac{1}{24}-\sum_{n=1}^{\infty}{\sigma(n)q^n} \quad \left(\sigma(n)\stackrel{\text{def}}{=}\sum_{d|n}{d}\right).
\end{equation}
It is easily seen that $\theta(G_{\eta}(\Gamma))=\theta(G(\Gamma))\cap \left(\bigoplus_{1<t|N}\Span_{\bQ}\{G_2-t(G_2|_2V_t)\}\right)$. \\

We now introduce the elements of the Hecke eigenbasis of $\cE_2(\Gamma)$. For $u\geq 1$ an integer, we let $\operatorname{P}(u)$ denotes the set of primitive Dirichlet characters of conductor $u$. For instance, $\operatorname{P}(1)=\{\textbf{1}_1\}$, where $\textbf{1}_1$ is the primitive Dirichlet character of conductor $1$, and $\operatorname{P}(2)$ is empty. If $u>1$ is such that $u^2|N$ and $\chi \in \operatorname{P}(u)$, the \textit{Eisenstein series attached to} $\chi$ is defined as
\begin{equation}
E_2^{\chi,\overline{\chi}}(\tau)\stackrel{\text{def}}{=}\sum_{n=1}^{\infty}{\sigma^{\chi,\overline{\chi}}(n)q^n} \quad \left(\sigma^{\chi,\overline{\chi}}(n)\stackrel{\text{def}}{=}\sum_{d|n}{\chi\left(\frac{n}{d}\right)\overline{\chi}(d)d}\right).
\end{equation}
Let $t$ be a positive divisor of $N/u^2$. When $\chi\neq \textbf{1}_1$, we define  $E_2^{\chi,\overline{\chi},t}$ as $(E_2^{\chi,\overline{\chi}}|V_t)$. If $\chi=\textbf{1}_1$, we define $E_2^{\chi,\overline{\chi},t}$ as $G_2-t(G_2|_2V_t)$. The notation above follows \cite{Diamond}.
The Hecke decomposition of $\cE(\Gamma)$ is given by 
\begin{equation}\label{eisenstein}
\cE_2(\Gamma)=\bigoplus_{u|s(N)}\bigoplus_{\chi \in \operatorname{P}(u)}\bigoplus_{1<tu^2|N} \Span_{\bC}\{E_2^{\chi,\bar{\chi},t}\}
\end{equation}
where $s(N)$ is the positive integer defined uniquely by the decomposition $N=ms(N)^2$ where $m$ is square-free (see \cite[Thm.~4.6.2]{Diamond}).

\section{Results and proofs}\label{sec:results}
From the explicit description \eqref{eisenstein} of the Eisenstein space together with Knopp and Mason's \cite[Thm.~2]{KnoMas}, we provide a simpler proof of \cite[Thm.~2]{WinKoh} and extend it to the case $s(N)=2$:
\begin{Theorem}\label{extension}
Let $\Gamma=\Gamma_0(N)$ and assume that $s(N)$ is $1$ or $2$. Then $G_c(\Gamma)=G_{\eta}(\Gamma)$. 
\end{Theorem}
\begin{proof}
As mentioned previously, $\operatorname{P}(1)=\{\textbf{1}_1\}$ and $\operatorname{P}(2)$ is empty. As $s(N)$ is either $1$ or $2$, \eqref{eisenstein} reduces to
\begin{equation}
\cE_2(\Gamma)=\bigoplus_{1<t|N} \Span_{\bC}\{E_2^{\textbf{1}_1,\textbf{1}_1,t}\}=\bigoplus_{1<t|N} \Span_{\bC}\{ G_2-t(G_2|_2V_t)\}. \nonumber
\end{equation}
Let $f$ be in $G_c(\Gamma)$ and write $\theta f=e+s$ where $e \in \cE_2(\Gamma)$ and $s \in \cS_2(\Gamma)$. There exists $h$ in $G_0(\Gamma)$ such that $\theta h=s$ and hence $g$ in $G_{\eta}(\Gamma)$ such that $\theta g=e$. In particular, $f/g=h$ belongs to $G_c(\Gamma)\cap G_0(\Gamma)$. By \eqref{liouville}, $f/g$ is a nonzero complex constant.
\end{proof}

For two integers $a$ and $b$, we denote $(a,b)$ their greatest common divisor. We set $\textbf{1}_{a|b}=1$ if $a$ divides $b$ and $\textbf{1}_{a|b}=0$ otherwise. For a Dirichlet character $\chi$, we denote $L(2,\chi)$ the value of the $L$-function of $\chi$ at $s=2$ and $g(\chi)$ the Gauss sum of $\chi$. 

\begin{Proposition}\label{expansion-at-cusp}
Let $N$ and $u>1$ be two positive integers such that $u^2|N$ and let $t$ be a positive divisor of $N/u^2$. Let $\chi\in \operatorname{P}(u)$. Let $\fc=a/c$ be a cusp of $\Gamma_0(N)$ written in lowest term ($c$ being nonzero). The constant term of $E_2^{\chi,\overline{\chi},t}$ at $\fc$ is given by 
\begin{equation}
-2\chi\left(-\frac{tac}{(c,t)^2 u}\right)\left(\frac{u}{2\pi}\right)^2  \frac{L(2,\chi^2)}{g(\chi)}\textbf{1}_{u|\frac{c}{(c,t)}}. \nonumber
\end{equation}
\end{Proposition}
\begin{proof}
Assume first that $t=1$ so that $E_2^{\chi,\overline{\chi},t}=E_2^{\chi,\overline{\chi}}$. By \cite[Section~4.6]{Diamond}, we have
\begin{equation}
-g(\chi)E_2^{\chi,\overline{\chi}}=\left(\frac{u}{2\pi}\right)^2\sum_{1\leq x,y,z\leq u}{\chi(x)\chi(y)G_2^{(xu,y+zu)}} \nonumber
\end{equation}
where, for $(c,d)$ in $\bZ^2$,
\begin{equation}
G_2^{(c,d)}(\tau)\stackrel{\text{def}}{=}\textbf{1}_{N|c}\zeta^{(d)}(2)-\left(\frac{2\pi}{N}\right)^2\sum_{n=1}^{\infty}{\sigma^{(c,d)}(n)q^{n/N}},
\end{equation}
with $\zeta^{(d)}(2)=\sum{n^{-2}}$ where $n$ runs over nonzero integers such that $n\equiv d \Mod{N}$, and $\sigma^{(c,d)}(n)=\sum{\operatorname{sgn}(m)m e^{2i\pi md/N}}$ where $m$ runs over positive and negative divisor of $n$ such that $n/m\equiv c\Mod{N}$.

Let $\gamma=\smallmat{a}{b}{c}{d}\in \SLZ$ such that it maps the cusp $\infty$ to $\fc$. We have
\begin{equation}
\left(G_2^{(xu,y+zu)}|_2\gamma\right)=G_2^{(xu,y+zu)\gamma}=\textbf{1}_{N|u(xa+zc)+yc}\zeta^{(u(xb+zd)+dy)}(2)+\cO(q) \nonumber
\end{equation}
where the first equality follows from \cite[Prop.~4.2.1]{Diamond}. In particular, the constant term is zero unless $u|c$. Then, without loss of generality, we assume that there is a positive integer $m$ such that $c=mu$. We are led to the computation of 
\begin{equation}\label{reduc-sum}
\sum_{1\leq x,y,z\leq u}{\chi(x)\chi(y)\textbf{1}_{N|u(xa+zc)+yc}\zeta^{(u(xb+zd)+dy)}(2)}.
\end{equation}
In the sum \eqref{reduc-sum} we operate the change of indexes $(x',z')=(x,z)\gamma$ which leads to 
\begin{equation}
\sum_{1\leq x',y,z'\leq u}{\chi(x'd-z'c)\chi(y)\textbf{1}_{N|ux'+yc}\zeta^{(uz'+dy)}(2)}, \nonumber
\end{equation}
and then to 
\begin{equation}\label{reduc-sum-two}
\chi(d)\sum_{1\leq x',y,z'\leq u}{\chi(x')\chi(y)\textbf{1}_{N|u(x'+ym)}\zeta^{(uz'+dy)}(2)}.
\end{equation}
The term $\textbf{1}_{N|u(x'+ym)}$ forces the relation $x'\equiv -ym\Mod{u}$ and \eqref{reduc-sum-two} restricts to 
\begin{equation}
\chi(-md)\sum_{1\leq y\leq u}{\chi(y)^2\left(\sum_{1\leq z'\leq u}{\zeta^{(uz'+dy)}(2)}\right)}. \nonumber
\end{equation}
As 
\begin{equation}
\sum_{z'\Mod{u}}{\zeta^{(uz'+dy)}(2)}=\sum_{\substack{n\in \bZ\setminus \{0\} \\ n\equiv dy\Mod{u}}}{\frac{1}{n^2}}, \nonumber
\end{equation}
and as $\{dy\Mod{u}~|~1\leq y\leq u\}$ is a complete set of representatives $\Mod{u}$, the sum \eqref{reduc-sum} equals $2\chi(-m)\overline{\chi}(d)L(2,\chi^2)$. As $m=c/u$ and $\overline{\chi}(d)=\chi(a)$, the constant term of $E^{\chi, \overline{\chi}}$ at $\fc$ is
\begin{equation}
-2\chi\left(-\frac{ac}{u}\right)\left(\frac{u}{2\pi}\right)^2\frac{L(2,\chi^2)}{g(\chi)}\textbf{1}_{u|c}. \nonumber
\end{equation}
The proposition for arbitrary $t$ follows from the next easy lemma:
\begin{Lemma}
Let $M$ and $t$ be two positive integers and $f$ be a modular form of integral weight $k$ for $\Gamma_0(M)$. The constant term of $(f|_kV_t)$ at a cusp $a/c$ of $\Gamma_0(Mt)$ is the constant term of $f$ at the cusp $ta/c$ of $\Gamma_0(M)$.
\end{Lemma}
It suffices to notice that $ta/c$ is $(t/(c,t))a/(c/(c,t))$ written in lowest term.
\end{proof}

Proposition \ref{expansion-at-cusp} motivates the computation of $L(2,\chi)$ when $\chi$ is an even Dirichlet character of conductor $p$. If $\chi$ is the trivial character, then $L(2,\chi)=\zeta(2)\prod_{p|N}{(1-p^{-2})}$. Otherwise, we record:
\begin{Lemma}\label{L-values}
Let $\chi$ be an even primitive Dirichlet character of conductor $p$. We have 
\begin{equation}
L(2,\chi)=\frac{1}{g(\overline{\chi})}\left(\frac{\pi}{p}\right)^2\left(\sum_{k=0}^{p}{k^2\overline{\chi}(k)}\right). \nonumber
\end{equation}
\end{Lemma}
\begin{proof}
The functional equation of $L(s,\chi)$ at $s=2$ (e.g. \cite[Thm.~12.11]{Apostol}) yields $(-2p)g(\overline{\chi})L(2,\chi)=(2\pi)^2 L(-1,\overline{\chi})$. If $\zeta(s,a)$ denotes the Hurwitz zeta function (see \cite[Chap.~12]{Apostol}), we have:
\begin{equation}
L(-1,\overline{\chi})=p\sum_{k=1}^p{\zeta\left(-1,\frac{k}{p}\right)\overline{\chi}(k)}=\frac{-1}{2p}\sum_{k=1}^{p}{k^2\overline{\chi}(k)} \nonumber
\end{equation}
using both \cite[Thm.~12.13]{Apostol} and the proof of \cite[Thm.~12.20]{Apostol}. We conclude.
\end{proof}

As a consequence of the previous computation, we otbain:
\begin{Theorem}\label{Gp}
Let $p$ be a prime number and $\Gamma=\Gamma_0(p^2)$. Let $f$ be an element of $G(\Gamma)$ having rational coefficients. Then $f$ belongs to $G_0(\Gamma)G_\eta(\Gamma)$.
\end{Theorem}

\begin{proof}
Let $f$ be an element in $G(\Gamma)$ having rational coefficients. Recall that $\theta$ induces a surjective group morphism from $G_{\eta}(\Gamma)$ to the intersection of $\bigoplus_{1<t|N} \operatorname{Span}_{\bC}\{E_2^{\textbf{1}_1,\textbf{1}_1,t}\}$ and $\theta(G(\Gamma))$. Hence by \eqref{eisenstein}, we need to show that the projection of $\theta f$ onto the direct summand $\bigoplus_{\chi \in \operatorname{P}(p)} \Span_{\bC}\{E_2^{\chi,\bar{\chi}}\}$ of $\cE_2(\Gamma)$ is zero. Let us write $g=\sum_{\chi\in \operatorname{P}(p)}{c_{\chi}E_2^{\chi,\overline{\chi}}}$ the latter projection, for coefficients $c_{\chi}\in \bC$. Since $\theta f$, $E_2^{\textbf{1}_1,\textbf{1}_1,p}$ and $E_2^{\textbf{1}_1,\textbf{1}_1,p^2}$ all have rational Fourier coefficients and constant terms at cusps, so is $g$:
\begin{equation}\label{pointes}
\sum_{\chi\in  \operatorname{P}(p)}{\chi(a)\frac{L(2,\chi^2)}{g(\chi)}c_{\chi}}\in \pi^2\bQ \quad (\forall a\in \bZ),
\end{equation}
\begin{equation} \label{coeff}
\sum_{\chi\in  \operatorname{P}(p)}{c_{\chi}\sigma^{\chi,\overline{\chi}}(n)}\in \bQ \quad (\forall n\geq 1). 
\end{equation}
Let $\psi$ be an element of $\operatorname{P}(p)$. Multiplying \eqref{pointes} by $\overline{\psi}(a)$ and summing over $a\in\{1,2,...,p-1\}$ gives
\begin{equation}
\frac{L(2,\psi^2)}{g(\psi)}c_{\psi}\in \pi^2\bQ(\psi) \nonumber
\end{equation}
where $\bQ(\psi)$ denotes the algebraic extension of $\bQ$ generated by the values of $\psi$. If $\psi$ is real then $\psi^2$ is trivial and $L(2,\psi^2)\in \pi^2\bQ$ implies that $c_{\psi}\in g(\psi)\bQ(\psi)$. Otherwise, $\psi^2$ is even and primitive and Lemma \ref{L-values}, we find $c_{\psi}\in g(\psi)g(\overline{\psi}^2)\bQ(\psi)$. Denoting $J(\overline{\psi},\overline{\psi})\in \bQ(\psi)$ the Jacobi sums associated to $(\overline{\psi},\overline{\psi})$, we have 
\begin{equation}
g(\psi)g(\overline{\psi}^2)=\psi(-1)p\frac{g(\overline{\psi})}{J(\overline{\psi},\overline{\psi})}\in g(\overline{\psi})\bQ(\psi).\nonumber
\end{equation}
It infers that regardless of $\psi$, $c_{\psi}\in g(\overline{\psi})\bQ(\psi)$.

Now, note that the arithmetic functions $\sigma^{\chi,\overline{\chi}}$ are linearly independent over $\bC$ (being the Fourier coefficients of linearly independent modular forms). Let $\zeta_{p(p-1)}$ be a primitive $p(p-1)$-root of unity. As the coefficients $c_{\chi}$ are in the cyclotomic field $\bQ(\zeta_{p(p-1)})$, Equation \eqref{coeff} implies that, for any $\rho$ in $\Gal(\bQ(\zeta_{p(p-1)})/\bQ)$ and $\chi$ in $\operatorname{P}(p)$, we have $\rho(c_{\chi})=c_{\rho(\chi)}$. Hence, $c_{\chi}$ belongs to $\bQ(\chi)$.

Then, if $c_{\psi}$ is non zero, it implies $g(\overline{\psi})\in \bQ(\psi)$ which is absurd: as $\psi$ has values in the $p-1$-roots of unity, there is an automorphism $\rho$ of $\Gal(\bQ(\zeta_{p(p-1)})/\bQ)$ leaving $\psi$ invariant but sending $e^{2i\pi/p}$ to $e^{2i\pi b/p}$ with $\psi(b)\neq 1$. As $\rho(g(\overline{\psi}))=\psi(b)g(\overline{\psi})$ this implies that $g(\overline{\psi})$ does not belong to $\bQ(\psi)$.
\end{proof}
\begin{Remark}
We believe that Theorem \ref{Gp} holds when $\Gamma_0(p^2)$ is replaced by any $\Gamma_0(N)$. But the present calculation gets sensitively harder and does not seem to be able to be completed without an additional argument. 
\end{Remark}
We are now in position to prove our main Theorem \ref{thmI}.
\begin{proof}[Proof of Theorem \ref{thmI}]
Clearly, \ref{item:i} implies \ref{item:ii} and \ref{item:ii} implies \ref{item:iii}. Also, \ref{item:i} implies \ref{item:iv} and \ref{item:iv} implies \ref{item:v}. If \ref{item:iii} is fulfilled, then by Theorem \ref{Gp}, there exists $g\in G_0(\Gamma)$ and $h\in G_{\eta}(\Gamma)$ normalized such that $f=gh$. It implies that $g=fh^{-1}$ has $q$-integral coefficients and is in $G_0(\Gamma)$. Then $g$ is constant by \cite[Thm~2]{KohMas}. Thus \ref{item:iii} implies \ref{item:i}. \\
If \ref{item:v} is fulfilled, we again write $f=gh$ with $g\in G_0(\Gamma)$ and $h\in G_{\eta}(\Gamma)$ by Theorem \ref{Gp}. Then, $\nu_{g}$ has absolute value one so $g$ is bounded on $\cH$. By \eqref{liouville}, it is constant. Hence, \ref{item:v} implies \ref{item:i}.
\end{proof}
If we do not assume the rationality of the coefficients, this is no longer true according Propositions \ref{holo} and \ref{classical} below.
\begin{Proposition}\label{holo}
Let $\chi$ be a primitive real Dirichlet character of conductor $u>1$. Let $\Gamma=\Gamma_0(u^2)$. The infinite product
\begin{equation}
\eta_{\chi}(\tau)\stackrel{\text{def}}{=}\prod_{n=1}^{\infty}{\prod_{\ell=0}^{u-1}{\left(1-e^{\frac{2i\pi \ell}{u}}q^n\right)^{\bar{\chi}(\ell n)}}}
\end{equation}
defines an element of $G(\Gamma)$ which is sent to $-(g(\bar{\chi})/2)E_2^{\chi,\bar{\chi}}$ under $\theta$.
\end{Proposition}
\begin{Remark}
The $\bar{\chi}$ in exponent in the definition of $\eta_{\chi}$ may seem unnecessary since we assume that $\chi$ is real. However, the proof infers that this choice would be natural if one defined $\eta_{\chi}$ for arbitrary $\chi$. 
\end{Remark}

\begin{proof}[Proof of Proposition \ref{holo}.]
Let $\log$ be the principal branch of the complex logarithm. The sum of holomorphic functions
\begin{equation}
\sum_{n=1}^{\infty}\sum_{\substack{\ell=1\\ (\ell,u)=1}}^u{\overline{\chi}(\ell n)\log\left(1-e^{\frac{2i\pi \ell}{u}}q^n\right)}\nonumber
\end{equation}
converges uniformly on compact subsets of $\cH$ to the $\log$ of $\eta_{\chi}$. This implies that $\eta_{\chi}$ is holomorphic on $\cH$ and does not vanish. We then have
\begin{equation}
\theta \eta_{\chi}=-\sum_{n=1}^{\infty}{\sum_{\substack{\ell=1\\ (\ell,u)=1}}^u{n\overline{\chi}(\ell n)\frac{e^{\frac{2i\pi \ell}{u}}q^n}{1-e^{\frac{2i\pi \ell}{u}}q^n}}}=-\sum_{n=1}^{\infty}{\sum_{m=1}^{\infty}\sum_{\substack{\ell=1\\ (\ell,u)=1}}^u{n\overline{\chi}(\ell n) e^{\frac{2i\pi m \ell}{u}}q^{nm}}}. \nonumber
\end{equation}
As the convergence of the latter is absolute, we can invert the summations. The claimed identity $\theta \eta_{\chi}=-(g(\bar{\chi})/2)E_2^{\chi,\bar{\chi}}$ follows as $\chi$ is primitive. \\
To deduce that $\eta_{\chi}$ belongs to $G(\Gamma)$, we need to show that the constant term of $\theta \eta_{\chi}$ at a cusp $\fc=a/c$ is in $(w_{\fc}c^2)^{-1}\bZ$. Having Proposition \ref{expansion-at-cusp} in mind, we can assume that $c=u$. The condition at $\fc$ becomes:
\begin{equation}\label{condition}
\chi(-a)\left(\frac{u}{2\pi}\right)^2 L(2,\chi^2)\in \frac{1}{u^2}\bZ. 
\end{equation}
As $\chi$ is real, $L(2,\chi^2)=\zeta(2)\prod_{p|u}{\left(1-p^{-2}\right)}$. As there is no primitive Dirichlet character modulo $2$, $u\geq 3$. It is then a simple exercise to verify that
\begin{equation}
\chi(-a)\frac{u^4}{(2\pi)^2} L(2,\chi^2)=\chi(-a)\frac{u^4}{24}\prod_{p|u}{\left(1-\frac{1}{p^2}\right)}\in \bZ. \nonumber
\end{equation}
This implies \eqref{condition} and concludes the proof.
\end{proof}
At this stage, it is absolutely not clear whether the multiplier system $\nu_{\eta_{\chi}}$, for $\chi$ a real primitive character of conductor $u>1$, is of finite order or not. We leave it as an open question. We are able to answer positively for the following functions: \\
Let $\chi_4$ be the non-trivial character (mod $4$). Proposition \ref{holo} implies that the following functions of $\tau\in \cH$ are in $G(\Gamma)$:
\begin{equation}
\eta_{\left(\frac{\cdot}{3}\right)}(\tau)=\prod_{n=1}^{\infty}{\left(\frac{1-e^{2i\pi/3}q^n}{1-e^{4i\pi/3}q^n}\right)^{\left(\frac{n}{3}\right)}} \quad \text{for}~\Gamma=\Gamma_0(9),
\end{equation}
\begin{equation}
\eta_{\chi_4}(\tau)=\prod_{n=1}^{\infty}{\left(\frac{1-iq^n}{1+iq^n}\right)^{\chi_4(n)}}\quad \text{for}~\Gamma=\Gamma_0(16),
\end{equation}
\begin{equation}
\eta_{\left(\frac{\cdot}{5}\right)}(\tau)=\prod_{n=1}^{\infty}{\left(\frac{1-\frac{1-\sqrt{5}}{2}q^n+q^{2n}}{1-\frac{1+\sqrt{5}}{2}q^n+q^{2n}}\right)^{\left(\frac{n}{5}\right)}}\quad \text{for}~\Gamma=\Gamma_0(25).
\end{equation}

\begin{Proposition}\label{classical}
The above elements of $G(\Gamma)$ are classical, i.e. belongs to $G_c(\Gamma)$.
\end{Proposition}
\begin{proof}
The genus of $X(u^2)$ is zero if, and only if, $u$ equals $3$, $4$ or $5$. A group $\Gamma$ produces a modular curve $X(\Gamma)$ of genus zero if, and only if, it is generated by parabolic and elliptic elements. In the latter case, any character $\nu:\Gamma\to \bC^{\times}$ is finite. This is then the case for $\nu_f$ when $f$ is $\eta_{\left(\frac{\cdot}{3}\right)}$, $\eta_{\chi_4}$ or $\eta_{\left(\frac{\cdot}{5}\right)}$.
\end{proof}

\begin{proof}[Proof of Theorem \ref{thmII}] 
It follows as a corollary of Propositions \ref{holo} and \ref{classical}.
\end{proof}

\end{document}